\documentclass{amsart}
\usepackage{amsmath}
\usepackage{amsthm}
\usepackage{amssymb}
\usepackage{amsbsy}
\usepackage{amsfonts}
\usepackage{amstext}
\usepackage{amscd}
\usepackage[dvips]{epsfig}

\numberwithin{equation}{section}
\theoremstyle{plain}
\newtheorem{thm}{Theorem}[section]
\newtheorem{prop}[thm]{Proposition}
\newtheorem{cor}[thm]{Corollary}
\newtheorem{lem}[thm]{Lemma}
\theoremstyle{definition}
\newtheorem{exa}[thm]{Example}

\newtheorem{rem}[thm]{Remark}
\newtheorem{defi}[thm]{Definition}
\DeclareMathOperator*{\real}{\mathbb{R}}

\DeclareMathOperator*{\comp}{\mathbb{C}}
\DeclareMathOperator*{\nat}{\mathbb{N}}
\newcommand{\pn}{\mathcal{P}(n)}
\newcommand{\ncpn}{\mathcal{NC}(n)}

\newcommand{\mpn}{\mathcal{M}(n)}

\begin{document}
\title{Joint cumulants for natural independence}

\author{Takahiro Hasebe}
\address{Graduate School of Science,  Kyoto University,  Kyoto 606-8502, Japan}
\thanks{TH is supported by Grant-in-Aid for JSPS Research Fellows.}
\email{hsb@kurims.kyoto-u.ac.jp}
\author{Hayato Saigo}
\address{Graduate School of Science,  Kyoto University,  Kyoto 606-8502, Japan}
\email{harmonia@kurims.kyoto-u.ac.jp}
\date{}

\maketitle

\begin{abstract}
Many kinds of independence have been defined in non-commutative probability theory. Natural independence is an important class of independence; this class consists of five independences (tensor, free, Boolean, monotone and anti-monotone ones). In the present paper, a unified treatment of joint cumulants is introduced for natural independence. The way we define joint cumulants enables us not only to find the monotone joint cumulants but also to give a new characterization of joint cumulants for other kinds of natural independence, i.e., tensor, free and Boolean independences. 

We also investigate relations between generating functions of moments and monotone cumulants. We find a natural extension of the Muraki formula, which describes the sum of monotone independent random variables, to the multivariate case.    
\end{abstract}

\section{Introduction}
Many kinds of independence are known in non-commutative probability theory. 
The most important example is the usual independence in probability theory, naturally extended to the non-commutative case. This is called tensor independence.  Free independence is another famous example \cite{V1, V2} and 
there are many researches on it (see \cite{V-D-N} for early results). After the appearance of free independence, Boolean \cite{S-W} and monotone independence \cite{Mur3} were found as other 
interesting examples of independence. To classify these independences, Speicher defined in \cite{Spe1} universal independence which satisfies some nice properties such as associativity of independence. After that, Sch\"{u}rmann and Ben Ghorbal formulated the universal independence in a categorical setting in \cite{B-S1}. In \cite{Mur0} Muraki defined quasi-universal 
independence which allows non-commutativity of independence by replacing partitions in the definition of universal independence by ordered partitions. 
Later Muraki introduced natural independence in \cite{Mur4} 
as a generalization of the paper \cite{B-S1}. He proved that there are only five kinds of natural independence: tensor, free, Boolean, monotone and anti-monotone independences. Since essential difference does not appear between monotone and anti-monotone independences for the purpose of this paper, we do not consider anti-monotone independence.

Let $(\mathcal{A}, \varphi)$ be an algebraic probability space, i.e., a pair of a unital $\ast$-algebra and a state on it. 
Let $\mathcal{A}_\lambda$ be $\ast$-subalgebras, where $\lambda \in \Lambda$ are indices.   
The above mentioned four independences are defined as rules to calculate moments $\varphi(X_1 \cdots X_n)$ for 
\[
X_i \in \mathcal{A}_{\lambda _i},~ \lambda_i \neq \lambda_{i+1},\: 1\leq i\leq n-1,\: n\geq 2. 
\]
\begin{defi}(1) Tensor independence: 
$\left\{\mathcal{A}_{\lambda} \right\}$ is tensor independent if
\[
\varphi(X_1\cdots X_n)= \prod_{\lambda \in \Lambda} \varphi\Big(\overrightarrow{\prod_{i; X_i \in \mathcal{A}_{\lambda}}}X_i \Big),   
\]
where $\overrightarrow{\prod}_{i \in V} X_i$ is the product of $X_i, ~i \in V$ in the same order as they appear in $X_1 \cdots X_n$. \\
(2) Free independence \cite{V1}: 
We assume all $\mathcal{A}_\lambda$ contain the unit of $\mathcal{A}$. 
$\left\{\mathcal{A}_{\lambda} \right\}$ is free independent if
\[
\varphi(X_1 \cdots X_n) = 0 
\]
holds whenever $\varphi (X_1)=\cdots=\varphi (X_n)=0$. \\
(3) Boolean independence \cite{S-W}: 
$\left\{\mathcal{A}_{\lambda} \right\}$ is Boolean independent if
\[
\varphi(X_1\cdots X_n)=\varphi(X_1) \cdots \varphi(X_n). 
\]
(4) Monotone independence \cite{Mur3}: 
We assume that $\Lambda $ is equipped with a linear order $<$. Then 
$\left\{\mathcal{A}_{\lambda} \right\}$ is monotone independent if
\[
\varphi (X_1\cdots X_i\cdots X_n)=\varphi (X_i)\varphi(X_1\cdots X_{i-1}X_{i+1}\cdots X_n)
\]
holds when $i$ satisfies $\lambda _{i-1} < \lambda _i$ and $\lambda _i > \lambda _{i+1}$ 
(one of the inequalities is eliminated when $i=1$ or $i=n$).
\end{defi} 
Independence for subsets $S_\lambda \subset \mathcal{A}$ is defined by taking 
the algebras $\mathcal{A}_\lambda$ generated by $S_\lambda$ (without the unit of $\mathcal{A}$ in the case of monotone or Boolean independence).

Many probabilistic notions have been introduced for each kind of independence. In particular, analogues of cumulants are a central topic in this field. In the usual probability theory, cumulants are extensively used in the study such as the correlation function of a stochastic process. When more than one random variables are concerned, cumulants for a single random variable are not adequate and their extension to the multivariate case is required. Cumulants for the multivariate case is called joint cumulants or sometimes multivariate cumulants. 
In free probability theory, Voiculescu introduced free cumulants in \cite{V1, V2} for a single random variable as an analogy of the cumulants in probability theory. Later Speicher defined free cumulants for the multivariate case \cite{Spe2}. Speicher also clarified that non-crossing partitions appear in the relation between moments and free cumulants. 
The reader is referred to \cite{N-S1} for further references. Boolean cumulants were introduced in \cite{S-W} in the single variable case and seemingly in \cite{Leh1} in the multivariate case. 

Lehner unified many kinds of cumulants in non-commutative 
probability theory in terms of Good's formula. 
A crucial idea was a very general notion of independence called an exchangeability system \cite{Leh1}. 
Monotone cumulants however cannot be defined in Lehner's approach. 
This is because monotone independence is non-commutative: 
if $X$ and $Y$ are monotone independent, then $Y$ and $X$ are not necessarily monotone independent. 
Therefore, the concept of ``mutual independence of random variables'' fails to hold.  
In spite of this, we found a way to define monotone cumulants uniquely for a single variable in \cite{H-S}. In the present paper, we generalize the method to define joint cumulants for monotone independence. 

For tensor, free and Boolean cumulants, the following properties are considered to be basic.   
\begin{itemize}
\item[(MK1)] Multilinearity: $K_n: \mathcal{A}^n \to \comp$ is multilinear.   
\item[(MK2)] Polynomiality: There exists a polynomial $P_n$  such that 
\[K_n(X_1, \cdots, X_n) = \varphi(X_1\cdots X_n) + P_n\big( \{ \varphi(X_{i_1}\cdots X_{i_p}) \}_{\substack{1 \leq p \leq n-1,\\ i_1 < \cdots < i_p}}\big).
\]  
\item[(MK3)] Vanishment: If $X_1, \cdots, X_n$ are divided into two independent parts, i.e., there exist nonempty, disjoint subsets $I$, $J \subset \{1,\cdots, n\}$ such that $I \cup J = \{1,\cdots, n \}$ and $\{X_i, i \in I\}$, $\{X_i, i \in J\}$  are independent, then $K_n(X_1, \cdots, X_n) = 0$.  
\end{itemize} 
Cumulants for a single variable can be defined from joint cumulants: $K_n(X) := K_n(X, \cdots, X)$. 
Clearly the additivity of cumulants for a single variable follows from the property (MK3): $K_n(X+Y) =K_n(X) + K_n(Y)$ if $X$ and $Y$ are independent. 

The additivity of monotone cumulants for a single variable does not hold because of the non-commutativity of monotone independence. Instead, we proved in \cite{H-S} that monotone cumulants for a single variable satisfy that $K_n^{M}(N.X_1) := K_n^{M}(X_1 + \cdots + X_N) = NK_n^{M}(X_1)$ holds if $X_1  \cdots, X_N$ are identically distributed and monotone independent.

The notion of a dot operation is important throughout this paper. This notion was used in the classical umbral calculus \cite{Rota}. 
Section \ref{GES} is devoted to the definition of the dot operation associated to each notion of independence. 

In Section \ref{gcumulant} we define joint cumulants for natural independence in a unified way along an idea similar to \cite{H-S}. 
The new notion here is monotone joint cumulants denoted as $K_n^M$. The property (MK3) however does not hold for the reason above. 
Alternatively, it is expected that (MK3) holds for identically distributed random variables in view of the single-variable case. 
This is, however, not the case; as we shall see later, $K_3^M(X, Y, X) \neq 0$ for monotone independent, identically distributed $X$ and $Y$. 
To solve this problem, we generalize the condition (MK3) in Section \ref{gcumulant}. We can prove the uniqueness 
of joint cumulants under the generalized condition.    

Then we prove the moment-cumulant formulae for natural independences in Section \ref{uniindep} and Section \ref{natindep}. The formulae for universal independences (tensor, free, Boolean) are known facts, but our proof relates the highest coefficients and the moment-cumulant formulae. This proof is however not applicable to the monotone case and monotone moment-cumulant formula is proved in a more direct way.   

In Section \ref{genfunc} we clarify the relation of generating functions for monotone independence. We need to introduce a parameter $t$ 
which arises naturally from the dot operation. This parameter can be understood to be a parameter of a formal convolution semigroup.

\section{Dot operation}\label{GES}
We used in \cite{H-S} the dot operation associated to a given notion of independence. This is also crucial in the definition 
of joint cumulants for natural independence, that is, tensor, free, Boolean and monotone ones. 
\begin{defi} \label{dot}
We fix a notion of independence among tensor, free, Boolean and monotone. 
Let $(\mathcal{A}, \varphi)$ be an algebraic probability space.  
We take copies $\{X ^{(j)}\}_{j \geq 1}$ in an algebraic probability space $(\widetilde{\mathcal{A}},\widetilde{\varphi})$ for every $X \in\mathcal{A}$
such that 

(1) $X \mapsto X^{(j)}$ is a $\ast$-homomorphism from $\mathcal{A}$ to $\widetilde{\mathcal{A}}$ for each $j \geq 1$; 

(2) $\widetilde{\varphi}(X_1^{(j)} X_2^{(j)} \cdots X_n^{(j)}) = \varphi(X_1X_2\cdots X_n)$ for any $X_i \in \mathcal{A}$, $j, n \geq 1$; 

(3) the subalgebras  
$\mathcal{A}^{(j)} := \{X ^{(j)}\}_{X \in \mathcal{A}}$, $j \geq 1$ are independent.  

\noindent
Then we define the dot operation $N.X$ 
by 
\[
N.X = X^{(1)} + \cdots + X^{(N)}
\]
for $X \in \mathcal{A}$ and a natural number $N \geq 0$. We understand that $0.X = 0$. 
Similarly we can iterate the dot operation more than once; for instance $N.(M.X)$ can be defined (in a suitable space). 
\end{defi}
\begin{rem}
(1) The notation $N.X$ is inspired from ``the classical umbral calculus'' \cite{Rota}. Indeed, this notion can be used to develop some kind of  
umbral calculus in the context of quantum probability. 

(2) In many cases, we denote $\widetilde{\varphi}$ by $\varphi$ for simplicity. 
\end{rem}

%%%%%%%%%%%%%%%%%%%%%%%%%%%%%%
We can explicitly construct the above copies as follows. Let $\star$ be any one of the natural products of states (tensor, free, Boolean and monotone) on the free product of algebras and 
$\Lambda:=\{(i_1,\cdots,i_n): i_j \in \nat ~(1 \leq j \leq n), n \in \nat\}$.  
For an algebraic probability space $(\mathcal{A}, \varphi)$, we prepare copies $\{(\mathcal{A}^{\lambda}, {\varphi}^{\lambda})\}_{\lambda \in \Lambda}$ of it, i.e., 
$(\mathcal{A}^{\lambda}, {\varphi}^{\lambda})=(\mathcal{A},\varphi)$ for any $\lambda \in \Lambda$. Let us define a free product of algebras $\widetilde{\mathcal{A}}:=\ast_{\lambda \in \Lambda}\mathcal{A}^{\lambda}$ 
and a natural product of states $\widetilde{\varphi}:=\star_{\lambda \in \Lambda }{\varphi}^{\lambda}$ on $\widetilde{\mathcal{A}}$. 
Let ${(\cdot)}^{\lambda}: \mathcal{A} \ni X \mapsto { X}^\lambda \in {\mathcal{A}}^\lambda \subset \widetilde{\mathcal{A}}$ be the embedding of 
$\mathcal{A}$ into $\widetilde{\mathcal{A}}$, where ${ X}^\lambda$ is equal to $X$ as an element of $\mathcal{A} = {\mathcal{A}}^\lambda$.  We denote by the same symbol ${(\cdot)}^{(i)}$ the map 
${\mathcal{A}}^{(i_1,\cdots,i_n)} \ni {X}^{(i_1,\cdots,i_n)} \mapsto {X}^{(i_1,\cdots,i_n,i)} \in {\mathcal{A}}^{(i_1,\cdots,i_n,i)} \subset \widetilde{\mathcal{A}}$, which can be extended to a $\ast$-homomorphism on 
$\widetilde{\mathcal{A}}$.  Then iteration of dot operations can be realized in this space. For instance, $N. (M.X)$ is defined as $\sum_{j=1}^N{\left(\sum_{i=1}^{M}{X}^{(i)}\right)}^{(j)}=\sum_{j=1}^N\sum_{i=1}^M X^{(i,j)}$.  
%%%%%%%%%%%%%%%%%%%%%%%%%%%%%%%%%%%%%%%%%% 

\begin{rem}
%(1) Lehner defined an exchangeability system to develop a general theory of joint cumulants in terms of Good's formula~\cite{Leh1}.  
%The above construction of $(\widetilde{\mathcal{A}},\widetilde{\varphi},)$ coincides with the exchangeability system for tensor, free and 
%Boolean independences. We however needed the larger space $(\mathcal{U}^{(\infty)}, \varphi^{(\infty)})$ to define the dot operation. \\
While tensor, free and Boolean independences provide exchangeability systems, monotone independence does not. However, we can extend an exchangeability system to include monotone independence. More precisely, an exchangeability system for an algebraic probability space $(\mathcal{A},\varphi)$ consists of copies $\{X^{(i)}\}_{i \geq 1}$ of random variables $X \in \mathcal{A}$ such that, for arbitrary random variables $X_1,\cdots,X_n \in \mathcal{A}$ and a sequence $(i_1,\cdots,i_n)$ of natural numbers, a joint moment $\varphi(X_1^{(i_1)} \cdots X_n^{(i_n)})$ is equal to $\varphi(X_1^{(\sigma(i_1))} \cdots X_n^{(\sigma(i_n))})$ under any permutation $\sigma$ of $\nat$. Let us consider a weaker invariance that the joint moment is invariant under any order-preserving permutation $\sigma$, i.e., a permutation $\sigma$ of $\nat$ such that $i < j$ implies $\sigma(i) < \sigma(j)$. 
Then the copies in Definition \ref{dot} satisfy this weaker invariance for monotone independence as well as for the other three independences. 
\end{rem}

\begin{prop} (Associativity of dot operation). We fix a notion of independence among the four. 
Then the dot operation satisfies that
\[
\varphi\big(N.(M.X_1) \cdots N.(M.X_n)\big) = \varphi\big((MN).X_1 \cdots (MN).X_n \big)
\]
for any $X_i \in \mathcal{A}$, $n \geq 1$. 
\end{prop}
\begin{proof}
$N.(M.X_i)$ is the sum 
\begin{equation}\label{eq001}
X_i^{(1,1)} + X_i^{(2,1)} + \cdots + X_i^{(M,1)} + X_i^{(1,2)} +\cdots + X_i^{(M,N)},  
\end{equation}
where $\{X_i^{(1,j)}\}_{i=1}^n, \cdots, \{X_i^{(M,j)}\}_{i=1}^n$ are independent for each $j$ and $\{X_i^{(1,j)} + X_i^{(2,j)} + \cdots + X_i^{(M,j)}\}_{i=1}^n$ ($j = 1, \cdots, N$) are independent. 
On the other hand, $(NM).X_i$ is the sum 
\begin{equation}\label{eq002}
X_i^{(1)} + \cdots + X_i^{(NM)}, 
\end{equation}
where $\{X_i^{(1)}\}_{i=1}^n, \cdots, \{X_i^{(NM)}\}_{i=1}^n$ are independent. Since natural independence is associative, 
the random variables in (\ref{eq002}) satisfy a stronger condition of independence than those in (\ref{eq001}).  By the way, the condition of 
independence in (\ref{eq001}) is enough to calculate the expectation only by sums and products of joint moments of $X_1, \cdots, X_n$.  Therefore, $\varphi\big(N.(M.X_1) \cdots N.(M.X_n)\big)$ must be equal to $\varphi\big((MN).X_1 \cdots (MN).X_n \big)$. 
\end{proof}

%\begin{exa} Calculation of joint moments of the monotone GES is 
%performed as follows. 
%\[
%\begin{split}
%\widetilde{\varphi} (X_1^{(2)} X_2^{(3)} X_3^{(2)}X_4^{(5)}) &=  \widetilde{\varphi}(X_1^{(2)} X_2^{(3)} X_3^{(2)}) \widetilde{\varphi} (X_4^{(5)}) \\
%                                                &= \widetilde{\varphi}(X_1^{(2)} X_3^{(2)}) \widetilde{\varphi}(X_2^{(3)}) \varphi (X_4) \\ 
%                                                   &=  \varphi(X_1 X_3) \varphi(X_2) \varphi (X_4).
%\end{split}
%\]
%\end{exa} 

\section{Generalized cumulants}\label{gcumulant}
The following properties are basic for joint cumulants in tensor, free and Boolean independences. 
\begin{itemize}
\item[(MK1)] Multilinearity: $K_n: \mathcal{A}^n \to \comp$ is multilinear.   
%%%% CAUTION : Reversed version, equivalence by induction in remark
\item[(MK2)] Polynomiality: There exists a polynomial $P_n$  such that 
\[K_n(X_1, \cdots, X_n) = \varphi(X_1\cdots X_n) + P_n\big( \{ \varphi(X_{i_1}\cdots X_{i_p}) \}_{\substack{1 \leq p \leq n-1,\\ i_1 < \cdots < i_p}}\big).
\]  
%%%%%
\item[(MK3)] Vanishment: If $X_1, \cdots, X_n$ are divided into two independent parts, i.e., there exist nonempty, disjoint subsets $I$, $J \subset \{1,\cdots, n\}$ such that $I \cup J = \{1,\cdots, n \}$ and $\{X_i, i \in I\}$, $\{X_i, i \in J\}$  are independent, then $K_n(X_1, \cdots, X_n) = 0$.  
\end{itemize}
%Uniqueness of such joint cumulants follows from Theorem \ref{unique} in which a stronger result is proved. 

Monotone cumulants do not satisfy (MK3), even if $X_i$ 's are identically distributed. 
For instance, $K_3 ^M(X, Y, X) = \frac{1}{2}(\varphi(X^2)\varphi(Y) - \varphi(X)\varphi(Y)\varphi(X))$ if $X$ and $Y$ are monotone independent 
(see Example \ref{exa112} in Section \ref{natindep}). 
Instead we consider the following property. 
\begin{itemize}
\item[(MK3')] Extensivity: $K_n(N.X_1, \cdots, N.X_n) = NK_n(X_1, \cdots, X_n)$. 
\end{itemize}
The terminology of extensivity is taken from the property of Boltzmann entropy. 

%\begin{rem}
%More generally, the following condition is enough to prove the uniqueness of cumulants. 
%\begin{itemize} 
%\item[(MK3'')\label{MK3''}] There exists a polynomial $Q_n$ without a constant or a linear term with respect to $N$ such that 
%\[
%K_n(N.X_1,  \cdots, N.X_n) = N K_n(X_1, \cdots, X_n) + Q_n(N, \{ \varphi(X_{i_1}\cdots X_{i_p}) \}_{\substack{1 \leq p \leq n-1,\\ i_1 < \cdots < i_p}}). 
%\]
%\end{itemize}
%There is no change in the proof and we do not consider this condition anymore in this paper. 
%\end{rem}
In the tensor, free and Boolean cases, it is well known that there exist cumulants which satisfy (MK1), (MK2) and (MK3), and hence generalized cumulants exist obviously. Here we discuss the uniqueness of generalized cumulants for all natural independences, including monotone independence.
%A key for the proof of uniqueness is the theorem below.
%%%%% CAUTION : by linearlity, N^2 appeear...
\begin{thm}\label{unique}
For any one of tensor, free, Boolean and monotone independences, joint cumulants satisfying (MK1), (MK2) and (MK3') are unique. 
\end{thm}
\begin{proof}
We fix a notion of independence. Let $\{K_n^{(1)}\}$ and $\{K_n^{(2)}\}$ be two families of 
cumulants with possibly different polynomials in the conditions (MK1), (MK2) and (MK3'). 
By the recursive use of (MK2), $\varphi(X_1 \cdots X_n)$ can be represented as a polynomial of $K_p^{(1)}$'s, and also as another polynomial of $K_p^{(2)}$'s: 
\begin{equation*}
\begin{split}
&\varphi(N.X_1 \cdots N.X_n) \\
&~~~= K_n ^{(1)}(X_1,\cdots, X_n) + Q^{(1)}_n(K_p^{(1)}(X_{i_1},\cdots, X_{i_p}):1 \leq p \leq n-1,~i_1 < \cdots < i_p) \\
&~~~= K_n ^{(2)}(X_1,\cdots, X_n) + Q^{(2)}_n(K_p^{(2)}(X_{i_1}, \cdots, X_{i_p}):1 \leq p \leq n-1,~i_1 < \cdots < i_p). \\
\end{split}
\end{equation*}
It follows from (MK1) that these polynomials $Q^{(1)}$ and $Q^{(2)}$ have no constant terms or linear terms with respect to $K_p^{(i)}$'s. Then $\varphi(N.X_1 \cdots N.X_n)$ has forms such as 
\begin{equation*}
\begin{split}
\varphi(N.X_1 \cdots N.X_n) &= NK_n ^{(1)}(X_1,\cdots, X_n) + \\ 
                              &~~N^2 \cdot (\text{a polynomial of $N$ and~} \{ K_p^{(1)}(X_{i_1},\cdots, X_{i_p}) \}_{\substack{1 \leq p \leq n-1,\\ i_1 < \cdots < i_p}} ) \\
            &= NK_n ^{(2)}(X_1,\cdots, X_n) + \\ 
                              &~~N^2 \cdot (\text{a polynomial of $N$ and~} \{ K_p^{(2)}(X_{i_1},\cdots, X_{i_p}) \}_{\substack{1 \leq p \leq n-1,\\ i_1 < \cdots < i_p}} )  
\end{split}
\end{equation*}
because both $K_p^{(1)}$'s and $K_p^{(2)}$'s satisfy (MK3'). 
The coefficients of $N$ in the above two lines must be the same. Therefore, $K_n^{(1)} = K_n ^{(2)}$ for any $n$. 
\end{proof}
%%%%%%%
The above theorem implies that generalized cumulants coincide with the usual cumulants in tensor, free and Boolean independences since (MK3') is weaker than (MK3). This is nothing but a new characterization of those cumulants. 

The existence of cumulants is not trivial. A key fact is the following. 

\begin{prop}\label{prop1}
For tensor, free, Boolean and monotone independence, $\varphi(N.X_1 \cdots N.X_n)$ is a polynomial of $N$ and 
$\varphi(X_{i_1} \cdots X_{i_k})$ $(1 \leq k \leq n, i_1 < \cdots < i_k)$ without a constant term with respect to N. 
\end{prop}
\begin{proof}
First we notice that there exists a polynomial $S_n$ (depending on the choice of independence) for any $n \geq 1$ 
such that if $\{X_i \}_{i = 1} ^n$ and $\{Y_j \}_{j=1}^n$ are independent, 
\begin{equation}\label{eq0}
\begin{split}
\varphi((X_1 + Y_1) \cdots (X_n + Y_n)) &= \varphi(X_1 \cdots X_n) + \varphi(Y_1 \cdots Y_n) \\ 
                      & ~~~~~~ + S_n \big( \{ \varphi(X_{i_1}\cdots X_{i_p}) \}_{\substack{1 \leq p \leq n-1,\\ i_1 < \cdots < i_p}}, \{\varphi(Y_{j_1}\cdots Y_{j_q}) \}_{\substack{1 \leq q \leq n-1,\\ j_1 < \cdots < j_q}} \big). 
\end{split}
\end{equation}
For each $i \in \{1,\cdots, n\}$, let $\{X_i ^{(j)}\}_{j \geq 1}$ be copies of $X_i$ appearing in Definition \ref{dot}. 
We prove the theorem by induction on $n$. The claim is obvious for $n = 1$ since the expectation is linear. 
We assume that the claim is the case for $n \leq k$. 
We replace $X_i$  and $Y_i$ in (\ref{eq0}) by $X_i ^{(1)}$ and  $X_i ^{(2)} + \cdots + X_i ^{(L+1)}$, respectively. Then one has 
\[
\begin{split}
&\varphi((L+1).X_1 \cdots (L+1).X_{k+1}) - \varphi(L.X_1 \cdots L.X_{k+1})  \\ 
                         &~~~~~~~~~~~~~~~~~~~~= \varphi(X_1 \cdots X_{k+1}) +  S_{k+1} \big( \{ \varphi(X_{i_1}\cdots X_{i_p}) \}_{\substack{1 \leq p \leq k,\\ i_1 < \cdots < i_p}}, \{\varphi(L.X_{j_1}\cdots L.X_{j_q}) \}_{\substack{1 \leq q \leq k,\\ j_1 < \cdots < j_q}} \big).  
\end{split}
\]
The right hand side is a polynomial of $L$ by assumption. 
Therefore, the sum 
\[
N \varphi(X_1 \cdots X_{k+1}) + \sum_{L = 0}^{N-1} S_{k+1} \big( \{ \varphi(X_{i_1}\cdots X_{i_p}) \}_{\substack{1 \leq p \leq k,\\ i_1 < \cdots < i_p}}, \{\varphi(L.X_{j_1}\cdots L.X_{j_q}) \}_{\substack{1 \leq q \leq k,\\ j_1 < \cdots < j_q}} \big)
\]
is also a polynomial of $N$ without a constant. 
\end{proof}

%On the other hand, the existence of generalized cumulants for monotone independence is not trivial. We will prove that the coefficient of $N$ in the proof above satisfies all conditions for our generalized cumulants.
\begin{defi}\label{defi1}
We define the $n$-th monotone (resp.\ tensor, free, Boolean) cumulant $K_n^M$ (resp.\ $K_n^T$, $K_n^F$, $K_n^B$) by the coefficient of $N$ 
in $\varphi(N.X_1 \cdots N.X_n)$ for monotone (resp.\ tensor, free, Boolean) independence. 
\end{defi}
It is easy to see from the proof of Proposition \ref{prop1} that the multilinearity (MK1) and polynomiality (MK2) hold. 
The extensivity (MK3') comes from the associative law of the dot operation as follows. 
\begin{prop}\label{MK3'}
The cumulants $K^M _n, K_n ^T, K_n ^F, K_n ^B$ satisfy the condition (MK3'). 
\end{prop}
\begin{proof}
The idea is the same as in \cite{H-S}. 
We recall that the dot operation is associative: 
\[
\varphi(M.(N.X_1)\cdots M.(N.X_n)) = \varphi((MN).X_1 \cdots (MN).X_n). 
\]
By definition, 
$\varphi(M.(N.X_1)\cdots M.(N.X_n))$ is of such a form as  
\[
MK_n(N.X_1, \cdots, N.X_n) + M^2 \cdot(\text{a polynomial of $M$ and~} \{ \varphi(N.X_{i_1}\cdots N.X_{i_p}) \}_{\substack{1 \leq p \leq n-1,\\ i_1 < \cdots < i_p}}) . 
\]
Also by definition $\varphi((MN).X_1 \cdots (MN).X_n)$ is of such a form as 
\[
MNK_n(X_1, \cdots, X_n) + M^2 N^2\cdot(\text{a polynomial of $MN$ and~} \{ \varphi(X_{i_1}\cdots X_{i_p}) \}_{\substack{1 \leq p \leq n-1,\\ i_1 < \cdots < i_p}}) . 
\] 
The coefficients of $M$ coincide, and hence, (MK3') holds. 
\end{proof}

We know that $K^T$, $K^F$ and $K^B$ are no other than the usual tensor, 
free and Boolean cumulants, respectively, because of Theorem \ref{unique}. 
Therefore, it is obvious that the property (MK3) holds. However, we can also prove (MK3) directly on the basis of Definition \ref{defi1} as follows. 
\begin{prop}\label{vanishing1}
The property (MK3) holds for tensor, free and Boolean independences.
\end{prop}
\begin{proof} We prove the claim for tensor independence; the other cases can be proved in the same way. Let $(\mathcal{A}_i, \varphi_i)$ be algebraic probability spaces for $i=1,2$ and $(\mathcal{A}_3, \varphi_3)$ be defined by 
$(\mathcal{A}_3, \varphi_3) = (\mathcal{A}_1 \ast \mathcal{A}_2, \varphi_1 \otimes \varphi_2)$.  Moreover, for $i=1,2,3$ let $(\widetilde{\mathcal{A}}_i, \widetilde{\varphi}_i, \{\iota_i^{(k)} \}_{k \geq 1})$ be the tensor exchangeability system constructed in \cite{Leh1}. Namely, let $\{(\mathcal{A}_i^{(k)}, \varphi_i^{(k)})\}_{k\geq 1}$ be copies of $(\mathcal{A}_i, \varphi_i)$ for each $i \in \{1,2,3\}$, 
 $\widetilde{\mathcal{A}}_i:=\ast_{k\geq 1} \mathcal{A}_i^{(k)}$, $\widetilde{\varphi}_i:=\otimes_{k \geq 1} \varphi_i^{(k)}$ and $\iota_i^{(k)}: \mathcal{A}_i \to \mathcal{A}_i^{(k)} \subset \widetilde{\mathcal{A}}_3$ be the natural inclusion. We shall prove that 
$\widetilde{\mathcal{A}}_1$ and $\widetilde{\mathcal{A}}_2$ are tensor independent in $(\widetilde{\mathcal{A}}_3, \widetilde{\varphi}_3)$. This follows from the equality of states 
\[
\widetilde{\varphi}_3 = \otimes_{k\geq 1} (\varphi_1^{(k)} \otimes \varphi_2^{(k)})= (\otimes_{k \geq 1} \varphi_1^{(k)})\otimes  (\otimes_{k \geq 1} \varphi_2^{(k)}) = \widetilde{\varphi}_1 \otimes \widetilde{\varphi}_2
\] 
under the natural isomorphism 
\[
\widetilde{\mathcal{A}}_3 = \ast_{k \geq 1}\big(\mathcal{A}^{(k)}_1 \ast \mathcal{A}_2^{(k)} \big) \cong \widetilde{\mathcal{A}}_1 \ast \widetilde{\mathcal{A}}_2. 
\] 
This is because the tensor product of states is commutative.  

Now we take $X_1, \cdots, X_n \in \mathcal{A}_1 \cup \mathcal{A}_2$ satisfying 
$I:=\{i; X_i \in \mathcal{A}_1\} \neq \emptyset$ and $J:=\{i; X_i \in \mathcal{A}_2\} \neq \emptyset$. 
Then, we have 
\[
\widetilde{\varphi}_3(N.X_1 \cdots N.X_n) = \widetilde{\varphi}_1\Big(\overrightarrow{\prod_{i \in I}} (N. X_i)\Big) \widetilde{\varphi}_2\Big(\overrightarrow{\prod_{j \in J}}(N.X_j)\Big),  
\] 
since the sets $\{N.X_i; i \in I \}$ and $\{N.X_i; i \in J \}$ are independent. 
The definition of cumulants and the property (MK3') imply that the left hand side contains the term $NK_n^T(X_1, \cdots, X_n)$ while the coefficient of $N$ in the right hand side is zero. Therefore, $K_n^T(X_1, \cdots, X_n)=0$. 
\end{proof}
\begin{cor}
For any one of tensor, free and Boolean independences, cumulants satisfying (MK1), (MK2) and (MK3) uniquely exist. 
\end{cor}

\section{New look at moment-cumulant formulae for universal independences}\label{uniindep}
Lehner proved in \cite{Leh1} the moment-cumulant formulae in a unified way for tensor, 
free and Boolean independence via Good's formula. 
Therefore, one may naturally expect that the moment-cumulant 
formulae can also be proved on the basis of Definition \ref{defi1}. 
In this section, the crucial concept is universal independence or a universal 
product introduced by Speicher in \cite{Spe1}. He proved that there are only three kind of universal independence, i.e., tensor, free and Boolean ones.  

We introduce preparatory notations and concepts. $\pi$ is said to be a partition of $\{1,\cdots, n \}$ if 
$\pi=\{V_1,\cdots,V_k\}$, where $V_i$ are non-empty, disjoint subsets of $\{1,\cdots, n \}$ and $\cup_{i=1}^k V_i = \{1,\cdots, n \}$. The number $k$ of elements of $\pi$ is denoted as $|\pi|$. A partition $\pi$ is said to be crossing if there are blocks $V,W \in \pi$ such that elements $a,c \in V$ and $b,d \in W$ exist satisfying $a<b<c<d$. $\pi$ is said to be
 non-crossing if it is not crossing. Moreover, a non-crossing partition $\pi$ is called an interval partition if there are natural numbers $0= m_1 < m_2 <\cdots < m_k < m_{k+1} =n$ such that $\pi=\{V_1,\cdots,V_k\}$, where $V_i=\{m_{i}+1, m_i+2,\cdots, m_{i+1}\}$ for $1 \leq i \leq k$. The sets of partitions, non-crossing partitions and interval partitions are respectively denoted as $\mathcal{P}(n)$, $\mathcal{NC}(n)$ and $\mathcal{I}(n)$.

A partial ordering can be defined on $\pn$. For partitions $\pi$ and $\sigma$, $\sigma \leq \pi$ means that for any block $V \in \sigma$, there exists a block $W \in \pi$ such that $V \subset W$. The partition consisting of one block $\{1,\cdots,n\}$ is larger than any other partition.  

For random variables $\{X_i \}_{i=1}^n$ and a subset $W=\{j_1, \cdots, j_k \}$ of $\{1,\cdots, n\}$ with $j_1 < \cdots <j_k$, let $X_W$ denote the product $\overrightarrow{\prod}_{i \in W} X_i = X_{j_1} \cdots X_{j_k}$.  
We use the same notation for multilinear functionals: for multilinear functionals $T_p: \mathcal{A}^p \to \comp$ ($1 \leq p \leq n$) and the subset $W$ above, we define $T_{k}(X_W):= T_k(X_{j_1},\cdots,X_{j_k})$. 
Moreover, for a partition $\pi = \{V_1, \cdots, V_{|\pi|}\}$ of $\{1,\cdots, n\}$, 
we define $T_\pi(X_1, \cdots, X_n)$ to be the product $T_{|V_1|}(X_{V_1}) \cdots T_{|V_{|\pi|}|}(X_{V_{|\pi|}})$.

Given a family $(\mathcal{A}_i, \varphi_i)$ and a partition  
$\pi =\{V_1, \cdots ,V_p \} \in \pn$, we denote $X_1 \cdots X_n \in \mathcal{A}_\pi$  when $X_i$ and $X_j$ are in 
the same $\mathcal{A}_k$ if $i$ and $j$ are in the same block of $\pi$. 
Consider a finer partition $\sigma =\{W_1, \cdots ,W_r \} 
\leq \pi$ and define $k(l)$ for $l=1,\cdots, r$ by $X_i \in \mathcal{A}_{k(l)}$ for $i \in W_l$.
In this case we put
\begin{equation}
\varphi^{\sigma}(X_1 \cdots X_n):=\varphi _{k(1)}(X_{W_1})\cdots\varphi _{k(r)}
(X_{W_r}).
\end{equation}

Let a product of states on (unital) algebras $\Big( (\mathcal{A}_1, \varphi_1), (\mathcal{A}_2, \varphi_2) \Big) 
\mapsto (\mathcal{A}_1 \ast \mathcal{A}_2, \varphi_1 \star \varphi_2)$ be given,   
where $\ast$ denotes the free product (with identification of units in the case of unital algebras). 
\begin{defi}
The product $\star$ is called a universal product if it satisfies the following properties. 
\begin{itemize}
\item[(1)] Associativity: For all pairs $(\mathcal{A}_1, \varphi_1)$, $(\mathcal{A}_2, \varphi_2)$ and $(\mathcal{A}_3, \varphi_3)$, 
\begin{equation}
\varphi_1 \star (\varphi_2 \star \varphi_3) = (\varphi_1 \star \varphi_2) \star \varphi_3 
\end{equation}
under the natural identification of $(\mathcal{A}_1 \ast\mathcal{A}_2) \ast \mathcal{A}_3$ with $\mathcal{A}_1\ast (\mathcal{A}_2 \ast\mathcal{A}_3)$. 
\item[(2)] Universal calculation rule for moments: There exist coefficients $c(\pi; \sigma) \in \comp$ depending on $\sigma \leq \pi \in \pn$ such that 
\begin{equation}
\varphi (X_1 \cdots X_n) = \sum_{\sigma \leq \pi}c(\pi; \sigma)\varphi^\sigma (X_1 \cdots X_n)
\end{equation} 
holds for any $\pi \in \pn$, $n \geq 1$ and any $X_1 \cdots X_n \in \mathcal{A}_\pi$. Here  
$\varphi$ stands for the product 
\[
\varphi = \varphi_{k_1}\star\varphi_{k_2} \star\cdots \star\varphi_{k_p} 
\]
if $X_1 X_2 \cdots X_n \in \ast^{p}_{i=1} \mathcal{A}_{k_i}$. 

\end{itemize}
The coefficients $c(\pi;\pi)$ are called the highest coefficients. 
\end{defi}

We give a new proof of the moment-cumulant formulae obtained in the literature.  
The proof below makes it clear how a partition structure appears in a moment-cumulant formula. 
The following lemma is a simple consequence of the condition (2) of a universal product and (MK2). 
\begin{lem}\label{effi}
Let $\star$ be a universal product, i.e., the tensor, free or Boolean product.  Then there exist $d(\pi) \in \comp$ for $\pi \in \pn$ such that 
\[
\varphi(X_1 \cdots X_n) = \sum_{\pi \in \pn}d(\pi)K_\pi  (X_1, \cdots, X_n). 
\]
\end{lem}
\begin{thm}\label{coe2}
Let $c(\pi; \sigma)$ be the universal coefficients for a given universal independence. 
Let $d(\pi)$ be as in Lemma \ref{effi}. Then $d(\pi) = c(\pi;\pi)$. 
\end{thm}
\begin{proof}
%We only prove the claim for tensor independence since the other cases can be proved in the same way. 
Let $\pi \in \pn$ and $X_1 \cdots X_n \in \mathcal{A}_\pi$. Then 
\[
\begin{split}
\varphi(N.X_1 \cdots N.X_n) 
&= \sum_{\sigma \leq \pi} c(\pi; \sigma)\varphi^\sigma (N.X_1 \cdots N.X_n) \\
&= c(\pi;\pi) N^{|\pi|}K_\pi(X_1, \cdots, X_n)  \\ 
&~~~~+ \text{a polynomial of $N$ with degree more than $|\pi|$.}
\end{split}
\]
On the other hand,  Lemma \ref{effi} implies that 
\[
\begin{split}
\varphi(N.X_1 \cdots N.X_n) 
&= \sum_{\sigma \in \pn} d(\sigma)K_\sigma (N.X_1, \cdots, N.X_n) \\
&= \sum_{\sigma \in \pn} d(\sigma)N^{|\sigma|} K_\sigma (X_1, \cdots, X_n).
\end{split}
\]
We used (MK3), or weaker, (MK3') in the second line. 
Then, by (MK3), which is stronger than (MK3'), $K_\sigma(X_1, \cdots, X_n) = 0$ unless $\sigma \leq \pi$. 
Therefore,    we have the form 
\[
\begin{split}
\varphi(N.X_1 \cdots N.X_n) &= d(\pi) N^{|\pi|}K_\pi(X_1, \cdots, X_n)  \\ &~~~~+ \text{a polynomial of $N$ with degree more than $|\pi|$.} 
\end{split}
\]
Since the coefficients of $N^{|\pi|}$ coincide, $d(\pi) = c(\pi;\pi)$. 
\end{proof}
We have used the vanishing property (MK3) of joint cumulants, not only (MK3'), 
for universal independence. Therefore, we cannot apply the above proof to monotone independence. We prove a moment-cumulant formula for monotone independence in the next section. 

The highest coefficients for tensor, free and Boolean products are known as follows.  
\begin{thm}\label{coe1}(R. Speicher \cite{Spe1}) 
The highest coefficients are given as follows.  
\begin{itemize}
\item[(1)] In the tensor case, $c(\pi;\pi) = 1$ for $\pi \in \pn$. 
\item[(2)] In the free case, $c(\pi;\pi) = 1$ for $\pi \in \ncpn$ and $c(\pi;\pi) = 0$ for $\pi \notin \ncpn$. 
 \item[(3)] In the Boolean case, $c(\pi;\pi) = 1$ for $\pi \in \mathcal{I}(n)$ and $c(\pi;\pi) = 0$ for $\pi \notin \mathcal{I}(n)$. 
\end{itemize}
\end{thm}
The above result, combined with Theorem \ref{coe2}, completes the unified proof for moment-cumulant formulae for universal products.     
Namely,  we obtain 
\begin{align}
&\varphi(X_1 \cdots X_n) = \sum_{\pi \in \pn}K_\pi^T  (X_1, \cdots, X_n), \\  
&\varphi(X_1 \cdots X_n) = \sum_{\pi \in \ncpn}K_\pi^F  (X_1, \cdots, X_n), \\  
&\varphi(X_1 \cdots X_n) = \sum_{\pi \in \mathcal{I}(n)}K_\pi^B  (X_1, \cdots, X_n). 
\end{align}

\section{The monotone moment-cumulant formula}\label{natindep}
We call a subset $V \subset \{1,\cdots, n\}$ a block of interval type if there exist $i,j$, $1 \leq i \leq n, 0 \leq j \leq n-i$ such that $V = \{i, \cdots, i + j\}$. We denote by $\mathit{IB}(n)$ the set of all blocks of interval type. 

Let $V$ be a subset of $\{1, \cdots, n \}$ written as $V = \{k_1, \cdots, k_{m}\}$ with $k_1 < \cdots < k_m$, $m = |V|$. We collect all $1 \leq i \leq m+1$ satisfying $k_{i-1} + 1 < k_{i}$, where $k_0 := 0$ and $k_{m+1}:=n+1$. We label them $i_1, \cdots, i_p$. Let $V_1, \cdots, V_{p}$ be blocks defined by $V_q:= \{k_{i_q -1} + 1, \cdots, k_{i_q} -1 \}$. 
The figure used in Theorem \ref{thm9} is helpful to understand the situation. 

Under the above notation, we can prove the following. 
\begin{prop}\label{eq2}
If $\{X_i \}_{i = 1} ^n$ and $\{Y_j \}_{j=1}^n$ are monotone independent, 
\begin{equation}
\begin{split}
\varphi((X_1 + Y_1) \cdots (X_n + Y_n)) &= \sum_{V \subset \{1,\cdots,n\} } \varphi(X_V) \prod_{j = 1}^{p} \varphi(Y_{V_j}). 
\end{split}
\end{equation}
\end{prop}
\begin{proof}
The subsets $V_j$ play roles of choosing positions of $Y_i$'s. Then the claim follows immediately. 
\end{proof}

Let us define a multilinear functional $\varphi_N(X_1,\cdots,X_n):= \varphi(N.X_1 \cdots N.X_n)$ for $n \in \nat$ and $N \in \nat$. Since this is a polynomial of $N$, we can replace $N \in \nat$ by $t \in \real$ and then obtain a multilinear functional $\varphi_t:\mathcal{A}^n \to \comp$ for $n \in \nat$ and $t \in \real$. As in Section \ref{uniindep}, let $\varphi_t(X_W)$ denote $\varphi_t(X_{j_1},\cdots,X_{j_k})$ for a subset $W=\{j_1, \cdots, j_k \}$ of $\nat$ with $j_1 < \cdots <j_k$. Then the following is immediate from Proposition \ref{eq2}. 
\begin{cor}\label{recurrence}
We have the following recurrent differential equations. \\
(1) $\frac{d}{dt}\varphi_t(X_1, \cdots, X_n) = \sum_{V \subset \{1,\cdots,n\}, V \neq \emptyset} K_{|V|}^M(X_V)\prod_{j = 1}^{p} \varphi_t(X_{V_j})$. \\
(2) $\frac{d}{dt}\varphi_t(X_1, \cdots, X_n) = \sum_{V \in IB(n)}K_{|V|}^M(X_V)\varphi_t(X_{V^c})$. 
\end{cor}
\begin{proof} We replace $X_i$ and $Y_i$ in Proposition \ref{eq2} by $N.X_i$ and $(N+M).X_i - N.X_i$ respectively. 
We notice that $\{N.X_i \}_{i = 1}^n$ and $\{(N+M).X_i - N.X_i \}_{i=1}^n$ are monotone independent and that $(N+M).X_i - N.X_i$ is identically distributed to $M.X_i$. We replace $N$ by $t$ and $M$ by $s$ and then the equality 
\[
\varphi_{t+s}(X_1, \cdots, X_n) = \sum_{V \subset \{1,\cdots,n\} } \varphi_t(X_V) \prod_{j = 1}^{p} \varphi_s(Y_{V_j}) 
\]
holds. The equations (1) and (2) follows from respectively the derivation $\frac{d}{dt}|_{t=0}$ and $\frac{d}{ds}|_{s=0}$. We note that the coefficient of $s$ appears only when $V^c \in IB(n)$ and therefore we obtain (2) by replacing $V^c$ by $V$.  
\end{proof}

Now we prove the moment-cumulant formula which generalizes the result for the single-variable case~\cite{H-S}. In addition to partitions, we need ordered partitions in this section. An ordered partition of $\{1,\cdots, n\}$ is a sequence $(V_1,\cdots, V_k)$, where $\{V_1,\cdots,V_k \}$ is a partition of $\{1,\cdots, n\}$. An ordered partition can be written as a pair $(\pi,\lambda)$, where $\pi$ is a partition and $\lambda$ is an ordering of the blocks. For blocks $V,W \in \pi$, we denote by $V >_\lambda W$ if $V$ is larger than $W$ under the order $\lambda$. Let $\mathcal{LP}(n)$ be the set of ordered partitions. 

For a non-crossing partition $\pi$, we introduce a partial order on $\pi$. For $V,W \in \pi$, $V \succ W$ means that there are $i,~j \in W$ such that $i < k < j$ for all $k \in V$. Visually $V \succ W$ means that $V$ lies in the inner side of $W$. We then define a subset $\mpn$ of $\mathcal{LP}(n)$ by 
\begin{equation}
\mpn:= \{(\pi, \lambda): \pi \in \ncpn,~ \text{if $V, ~W \in \pi$ satisfy $V \succ W$,  then $V >_{\lambda} W$} \}.  
\end{equation}
An element of $\mathcal{M}(n)$ is called a monotone partition. 
The set of monotone partitions was first introduced by Muraki~\cite{Mur0} to classify natural independence. 
\begin{thm}
The moment-cumulant formula is expressed as 
\[
\varphi(X_1 \cdots X_n) = \sum_{(\pi, \lambda) \in \mpn} \frac{1}{|\pi|!}K_\pi ^M (X_1, \cdots, X_n).   
\]
\end{thm}
\begin{proof}
We prove this by induction on $n$. Assume that 
\[
\varphi_t(X_1, \cdots, X_k) = \sum_{(\pi, \lambda) \in \mathcal{M}(k)} \frac{t^{|\pi|}}{|\pi|!}K_\pi^M  (X_1, \cdots, X_k) 
\]
holds for $t \in \real$ and $k \leq n$. 
We recall that an element $(\pi,\lambda) \in \mathcal{M}(n)$ can be expressed as a sequence $(V_1, \cdots, V_{|\pi|})$. We can use a discussion similar to \cite{Has3, H-S}. A prototype of this discussion is in \cite{Saigo}. Let $IB(k, m)$ be the subset of $IB(k)$ defined by $\{V \in IB(k); |V|=m \}$. Let $1_k$ be the partition of $\mathcal{P}(k)$ consisting of one block.  There is a bijection $f: \mathcal{M}(n+1) \to \Big(\bigcup_{k = 1} ^{n} \mathcal{M}(n+1 - k)\times IB(n+1, k) \Big) \cup \{1_{n+1} \}$ defined by 
\[
f: (V_1, \cdots, V_{|\pi|}) \mapsto ((V_1, \cdots, V_{|\pi| - 1}), V_{|\pi|}). 
\] 
Therefore, the sum $\sum_{(\pi, \lambda) \in \mpn}$ 
can be replaced by $\sum_{V \in IB(n+1)} \sum_{(\sigma, \mu) \in \mathcal{M}(n+1 - |V|)}$ and we have 
\[
\begin{split}
\sum_{(\pi, \lambda) \in \mathcal{M}(n+1)} \frac{t^{|\pi|}}{|\pi|!}K_\pi ^M (X_1, \cdots, X_n) &= \sum_{V \in IB(n+1)} \sum_{(\sigma, \mu) \in \mathcal{M}(n +1- |V|)}\frac{t^{|\sigma|+1}}{(|\sigma|+1)!}K_\sigma^M (X_{V^c})K_{|V|}^M(X_V)  \\
                                                 &=  \sum_{V \in IB(n+1)} \int_0 ^t ds \sum_{(\sigma, \mu) \in \mathcal{M}(n+1 - |V|)}\frac{s^{|\sigma|}}{|\sigma|!}K_\sigma ^M(X_{V^c})K_{|V|}^M(X_V)  \\ 
                            &=  \sum_{V \in IB(n+1)} \int_0 ^t ds \varphi_s (X_{V^c}) K_{|V|}^M(X_V)  \\ 
                            &=   \int_0 ^t \frac{d}{ds}\varphi_s(X_1, \cdots, X_{n+1})ds \\ 
                            &=  \varphi_t(X_1, \cdots, X_{n+1}). 
\end{split}
\]
We used assumption of induction in the third line and Corollary \ref{recurrence} (2) in the fourth line. 
The claim follows from the case $t = 1$. 
\end{proof}

\begin{exa}\label{exa112}
We show the monotone cumulants up to the forth order. 
\[
\begin{split}
&K_1 ^M(X_1) = \varphi(X_1),~~~~ K_2 ^M(X_1, X_2) = \varphi(X_1X_2) - \varphi(X_1)\varphi(X_2), \\
&K_3 ^M(X_1, X_2, X_3) = \varphi(X_1X_2X_3) - \varphi(X_1X_2)\varphi(X_3) - \varphi(X_1)\varphi(X_2X_3) -\frac{1}{2}\varphi(X_1X_3)\varphi(X_2) \\
&~~~~~~~~~~~~~~~~~~~~~~~~~~~~~+ \frac{3}{2}\varphi(X_1)\varphi(X_2)\varphi(X_3), \\
&K_4 ^M(X_1,X_2,X_3,X_4) = \varphi(X_1X_2X_3X_4)- \varphi(X_1X_2X_3)\varphi(X_4)- \frac{1}{2}\varphi(X_1X_3X_4)\varphi(X_2) \\&~~~~~~~~~~~~~~~~~~~~~~~~~ - \frac{1}{2}\varphi(X_1X_2X_4)\varphi(X_3) 
  -\varphi(X_1)\varphi(X_2X_3X_4)  -\varphi(X_1X_2)\varphi(X_3X_4) \\&~~~~~~~~~~~~~~~~~~~~~~~~~ -\frac{1}{2}\varphi(X_1X_4)\varphi(X_2X_3) +\frac{3}{2}\varphi(X_1X_2)\varphi(X_3)\varphi(X_4)+\frac{2}{3}\varphi(X_1X_4)\varphi(X_2)\varphi(X_3) \\&~~~~~~~~~~~~~~~~~~~~~~~~~ + \frac{3}{2}\varphi(X_1)\varphi(X_2)\varphi(X_3X_4) +\frac{1}{2}\varphi(X_1)\varphi(X_2X_4)\varphi(X_3) + \frac{3}{2}\varphi(X_1)\varphi(X_2X_3)\varphi(X_4) \\&~~~~~~~~~~~~~~~~~~~~~~~~~+ \frac{1}{2}\varphi(X_1X_3)\varphi(X_2)\varphi(X_4) -\frac{8}{3} \varphi(X_1)\varphi(X_2)\varphi(X_3)\varphi(X_4). 
\end{split} 
\]
\end{exa}

\section{Generating functions}\label{genfunc}
Let $\comp[[z_1, \cdots, z_r]]$ be the ring of formal power series of non-commutative generators $z_1, \cdots, z_r$. An element $P(z_1, \cdots, z_r)$ in  $\comp[[z_1, \cdots, z_r]]$ 
can be expressed as 
\[
P(z_1, \cdots, z_r) =  p_{\emptyset} + \sum_{n=1}^\infty \sum_{i_1, \cdots, i_n = 1} ^r p_{i_1, \cdots, i_n}z_{i_1}\cdots z_{i_n}. 
\]
We define a generating function of the joint moments of $X=(X_1, \cdots, X_r)$ by 
\[
M_X (z_1, \cdots, z_r) := 1 + \sum_{n=1}^\infty \sum_{i_1, \cdots, i_n = 1} ^r \varphi(X_{i_1}\cdots X_{i_n})z_{i_1}\cdots z_{i_n} \in \comp[[z_1, \cdots, z_r]]. 
\]
First we show the following ``multivariate Muraki formula'' for generating functions. 
\begin{thm}\label{thm9} For any $X=(X_1, \cdots, X_r)$ and $Y=(Y_1, \cdots, Y_r)$ with $\{X_i \}_{i=1}^{r}$ and $\{Y_j\}_{j=1}^{r}$ monotone independent, 
\[M_{X+Y} (z_1,\cdots, z_r)=M_Y(z_1,\cdots,z_r)M_X (z_1 M_Y(z_1, \cdots, z_r),\cdots, z_r M_Y(z_1, \cdots, z_r) ).  
\] 
\end{thm}
\begin{proof} 
 For a fixed sequence $(i_1,\cdots,i_n)$, $1 \leq i_1, \cdots, i_n \leq r$, let us compare the coefficient of $z_{i_1} \cdots z_{i_n}$ in the both hands sides. In the left hand side, it was calculated in Proposition \ref{eq2}. The right hand side is expanded as 
\[
\begin{split}
&M_YM_X(z_1M_Y, \cdots, z_rM_Y) \\
&~~~~~= \sum_{k=0}^\infty \sum_{j_1, \cdots, j_k = 1} ^r \varphi(X_{j_1}\cdots X_{j_k})M_Y z_{j_1}M_Y z_{j_2}M_Y\cdots z_{j_k}M_Y, \\
\end{split}
\]
where the summation is understood to be $M_Y$ for $k=0$. The question is when the term $z_{i_1} \cdots z_{i_n}$ appears in $M_Y z_{j_1}M_Y z_{j_2}M_Y\cdots z_{j_k}M_Y$. This happens if and only if  the sequence $(j_1,\cdots, j_k)$ is a subsequence of $(i_1,\cdots, i_n)$. In this case, we can interpolate $(j_1,\cdots,j_k)$ to recover the whole sequence $(i_1,\cdots,i_n)$, by choosing unique terms from $M_Y$'s appearing in $M_Y z_{j_1}M_Y z_{j_2}M_Y\cdots z_{j_k}M_Y$.  In terms of a partition of a set $\{i_1,\cdots,i_n\}$, $(j_1, \cdots, j_k)$ can be described by a block $V$ and then the other blocks $(V_i)$ as in Fig.~\ref{fig1} interpolate $(j_1, \cdots, j_k)$. From Proposition \ref{eq2}, the coefficients of the both hands sides coincide. 
\begin{figure}[!htb]
\centering
\includegraphics[width=6.5cm,clip]{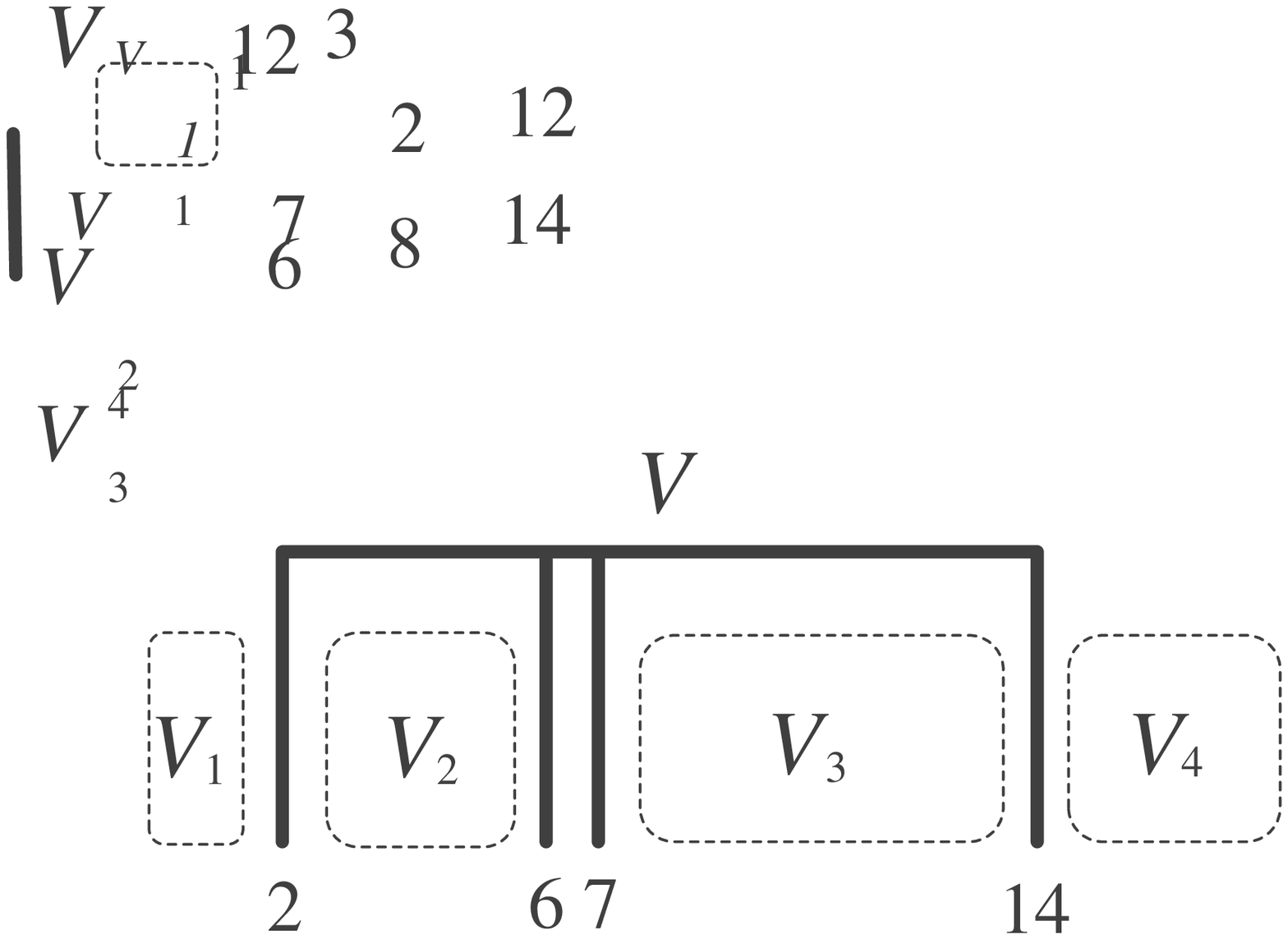}
\caption{This figure corresponds to the expectation $\varphi(Y_1)\varphi(Y_3Y_4Y_5)\varphi(Y_8 \cdots Y_{13})\varphi(Y_{15}Y_{16}Y_{17})\varphi(X_2 X_6X_7X_{14})$. The blocks $V_1,V_2,V_3,V_4$ are defined by $V_1=\{1\}, V_2=\{3,4,5\}, V_3=\{8,9,10,11,12,13\}$ and $V_4=\{15,16,17\}$.}
\label{fig1}
\end{figure}
\end{proof}

A generating function of the monotone cumulants of $X=(X_1, \cdots, X_r)$ is defined by  
\[
K_{X}^M(z_1, \cdots, z_r) := \sum_{n=1}^\infty \sum_{i_1, \cdots, i_n = 1} ^r K_n^M(X_{i_1},\cdots, X_{i_n})z_{i_1}\cdots z_{i_n} \in \comp[[z_1, \cdots, z_r]].\]
We denote by $M_{X}(t; z_1, \cdots, z_r)$ the generating function of the joint moments for the multilinear functionals $\varphi_t(X_1,\cdots,X_n)$. We also denote $M_{X}(z_1, \cdots, z_r)$ simply by $M_X$; $M_{X}(t; z_1, \cdots, z_r)$ by $M_X(t)$; $K^M_{X}(z_1, \cdots, z_r)$ by $K_X^M$. 
An important property is that $\frac{\partial M_X(t)}{\partial t}|_{t=0} = K_X^M$ holds. 

For random variable $X=(X_1,\cdots,X_r)$, let $\mu_{X,i}(z_1, \cdots, z_r) := z_i M_X(z_1, \cdots, x_r)$, $\mu_{X,i}(t):=z_iM_X(t)$ and $\kappa_{X,i}(z_1, \cdots, z_r) := z_i K_X^M(z_1, \cdots, z_r)$. We also introduce vectors $\mu_X:=(\mu_{X,1}, \cdots \mu_{X,r})$, $\mu_X(t):=(\mu_{X,1}(t),\cdots,\mu_{X,r}(t))$ and $\kappa_X:=(\kappa_{X,1}, \cdots \kappa_{X,r})$. One can see that every component of a vector has the same information. Therefore, one component is sufficient to understand the whole information on joint moments or cumulants. However, these vectors are useful to formulate a ``multivariate Muraki's formula''. 
\begin{cor}\label{cor8}
For any $X=(X_1, \cdots, X_r)$ and $Y=(Y_1, \cdots, Y_r)$ where $\{X_i \}_{i=1}^{r}$ and $\{Y_i\}_{i=1}^{r}$ are monotone independent, 
\[
\mu_{X+Y} = \mu_X \circ \mu_Y. 
\]
\end{cor}
Using this, we can derive a relation between a flow and a vector field. 
\begin{thm}\label{diffeq0} 
The following equalities hold. 
\begin{itemize}
\item[(1)] $\mu_X(t + s) = \mu_X(t) \circ \mu_X(s)$. 
\item[(2)] $\frac{\partial M_X(t)}{\partial t} = M_X(t) K_X^M(z_1 M_X(t), \cdots, z_r M_X(t))$, or equivalently, $\frac{\partial \mu_{X}(t)}{\partial t} = \kappa_X (\mu_X(t))$. 
\end{itemize} 
\end{thm}
\begin{proof}
(1) is immediate from Corollary \ref{cor8}: one just has to replace $X$ by $X^{(1)}+\cdots +X^{(M)}$ and $Y$ by $X^{(M+1)} + \cdots +X^{(M+N)}$. Then (1) is true as formal power series, where coefficients are polynomials regarding $M$ and $N$. Then we can extend $N$ and $M$ to real numbers $t$ and $s$, respectively. (2) follows from the derivative $\frac{d}{dt}|_0$. 
\end{proof}
It is worthy to compare Theorem \ref{diffeq0}(2) with the relation in free probability. 
Let $R_X(z_1, \cdots, z_r)$ be the generating function of free cumulants 
\[
R_X(z_1, \cdots, z_r) := \sum_{n=1}^\infty \sum_{i_1, \cdots, i_n = 1} ^r K_n^F(X_{i_1},\cdots, X_{i_n})z_{i_1}\cdots z_{i_n} \in \comp[[z_1, \cdots, z_r]]. 
\]
Then it is known that 
\begin{equation}\label{fcum}
M_X -1 = R_X(z_1 M_X, \cdots, z_r M_X).
\end{equation}
The reader is referred to Corollary 16.16 in \cite{N-S1}. %The right hand side is of a form similar to that in Theorem \ref{diffeq0} (1). Indeed, 
The above relation can also be expressed as $M_X -1 = R_X \circ \mu_X$ which is similar to the differential equation in Theorem \ref{diffeq0}(2).  
\begin{rem}\label{motivation}
 In the previous paper \cite{H-S}, we did not mention the relation between generating functions and cumulants. 
Now we explain the relation in detail. 
The differential equation becomes $\frac{\partial }{\partial t}M_X(t; z) = M_X(t; z) K_X^M(zM_X(t; z))$ in the one variable case. 
If we use $A_X(z):= -zK_X^M(\frac{1}{z})$ and the reciprocal Cauchy transform $H_X(t;z) = \frac{z}{M(t;\frac{1}{z})}$, the differential equation becomes 
\begin{equation}\label{diffeq1}
\frac{\partial }{\partial t} H_X (t;z) = A_X(H_X(t;z)). 
\end{equation}
This is the basic relation of a monotone convolution semigroup, first obtained in \cite{Mur3}.
Actually, a motivation of the paper \cite{H-S} was the observation that the coefficients of $A_X(z)$ had nice properties as cumulants. For instance, the arcsine law 
with mean $0$ and variance $1$ is characterized by $A_X(z) = -\frac{1}{z}$, or equivalently, $K^M_1(X) = 0$, $K^M_2(X) = 1$, $K^M_n(X) = 0$ for $n \geq 3$. Therefore, 
the problem was how to define cumulants for all probability measures. We can say that we defined monotone cumulants so that (\ref{diffeq1}) holds. In a recent paper \cite{Has3}, another way is presented to define monotone cumulants and their generalization on the basis of the differential equation (\ref{diffeq1}). However, it is difficult to generalize the method in \cite{Has3} to the multivariate case. In this sense, the present method has advantage. Theorem \ref{diffeq0} extends (\ref{diffeq1}) to the multivariate case. 

As is explained in the above, $t$ means a parameter of a ``formal'' convolution semigroup. Let us focus on this point more. Let $X$ be bounded and self-adjoint for simplicity. Then $M_X(t;z)$ may not be a moment generating function of a probability measure for general $t \geq 0$ and $X$. More precisely, $M_X(t;z)$ becomes a moment generating function of a probability measure for any $t \geq 0$ if and only if the probability distribution of $X$ is monotone infinitely divisible.  

The reader might wonder if there is a relation between the moment and cumulant generating functions without the use of $t$. For instance, one does not need the parameter $t$ in free probability theory \cite{V2}. In this case the cumulant generating function $K_X$ is called an $R$-transform and is denote by $R_X$. The basic relation is given by 
\[
M_X(z) = 1 + R_X(zM_X(z)). 
\] 
Therefore, $R_X$ can be expressed by using the inverse function of $zM_X(z)$. 
However, such a relation does not exist for monotone cumulants because of the difficulty of the correspondence between a holomorphic map and its vector field 
\cite{Bel, Cow, Has2}. 

In spite of the above, we can also understand this difficulty in a positive way since the use of the parameter $t$ indicates a new insight into relationship between independence and differential equations. 
%For instance, a complex Burgers equation was derived in free probability theory for a free convolution semigroup \cite{V2}. 
%In spite of this, we can understand the strong relation between cumulants and infinite divisibility. 
%In spite of the above, one can see a similarity between free probability and monotone probability:  if the time derivative $\frac{\partial}{\partial t}$ is erased and if $M_t$ is replaced by $M$, $\frac{\partial M_t}{\partial t} = M_t K(z_1M_t, \cdots, z_rM_t)$ becomes the relation in free probability \cite{}.  
\end{rem}

\section*{Acknowledgements} 
The authors thank Professor Izumi Ojima, Mr.\ Ryo Harada, 
Mr.\ Hiroshi Ando, Mr.\ Kazuya Okamura for discussions on the notion of independence. TH thanks Dr.\ Jiun-Chau Wang for a question on monotone cumulants and generating functions, which was a reason for our having written Remark \ref{motivation}. TH is supported by JSPS (KAKENHI 21-5106).


\begin{thebibliography}{122}
\bibitem{Bel} S.T.\ Belinschi, Complex analysis methods in noncommutative probability, PhD thesis, Indiana University, 2005, arXiv:math/0602343v1. 
     \bibitem{Cow} C.C.\ Cowen, Iteration and the  solution of functional equations for functions analytic in the unit disk, Trans.\ Amer.\ Math.\ Soc.\ \textbf{265}, No.\ 1 (1981), 69--95. 
\bibitem{B-S1} A.\ Ben Ghorbal and M.\ Sch\"{u}rmann, Non-commutative notions of stochastic independence, Math.\ Proc.\ Comb.\ Phil.\ Soc.\ \textbf{133} (2002), 531--561.   
    \bibitem{Has2} T.\ Hasebe, Monotone convolution and monotone infinite divisibility from complex analytic viewpoint, Infin.\ Dimens.\ Anal.\ Quantum Probab.\ Relat.\ Top.\ \textbf{13}, No.\ 1 (2010), 111--131.   
 \bibitem{Has3} T.\ Hasebe, Conditionally monotone independence I: Independence, additive convolutions and related convolutions, Infin.\ Dimens.\ Anal.\ Quantum Probab.\ Relat.\ Top., to appear. arXiv:0907.5473v4.    
  \bibitem{H-S} T.\ Hasebe and H.\ Saigo, The monotone cumulants, to appear in Ann.\ Inst.\ Henri Poincar\'{e} Probab.\ Stat., arXiv:0907.4896v3.   
   \bibitem{Leh1} F.\ Lehner, Cumulants in noncommutative probability theory I, Math.\ Z.\ \textbf{248} (2004), 67--100.   
  \bibitem{Mur3} N.\ Muraki, Monotonic convolution and monotonic L\'{e}vy-Hin\v{c}in formula, preprint, 2000. 
  \bibitem{Mur0} N.\ Muraki, The five independences as quasi-universal products, Infin.\ Dimens.\ Anal.\ Quantum Probab.\ Relat.\ Top.\ \textbf{5}, No.\ 1 (2002), 113--134.
     \bibitem{Mur4} N.\ Muraki, The five independences as natural products, Infin.\ Dimens.\ Anal.\ Quantum Probab.\ Relat.\ Top.\ \textbf{6}, No.\ 3 (2003), 337--371. 
     \bibitem{N-S1} A.\ Nica and R.\ Speicher, \textit{Lectures on the Combinatorics of Free Probability}, London Math.\ Soc.\ Lecture Note Series, vol.\ 335, Cambridge Univ.\ Press, 2006.   
     \bibitem{Rota} G.-C.\ Rota and B.D.\ Taylor, The classical umbral calculus, SIAM J.\ Math.\ Anal.\ \textbf{25} (1994), 694--711. 
   \bibitem{Saigo} H.\ Saigo, A simple proof for monotone CLT, Infin.\ Dimens.\ Anal.\ Quantum Probab.\ Relat.\ Top.\  \textbf{13}, No.\ 2  (2010), 339--343. 
   \bibitem{Spe2} R.\ Speicher, Multiplicative functions on the lattice of non-crossing partitions and free convolution, 
    Math.\ Ann.\ \textbf{298} (1994), 611--628.
      \bibitem{Spe1} R.\ Speicher, On universal products,  in Free Probability Theory, ed.\ D.\ Voiculescu, Fields Inst.\ Commun.\ \textbf{12} (Amer.\ Math.\ Soc., 1997), 257--266.  
    \bibitem{S-W} R.\ Speicher and R.\ Woroudi, Boolean convolution, in Free Probability Theory, ed.\ D.\ Voiculescu, Fields Inst.\ Commun.\ \textbf{12} (Amer.\ Math.\ Soc., 1997), 267--280. 
\bibitem{V1} D.\ Voiculescu, Symmetries of some reduced free product algebras, Operator algebras and their connections with topology and 
    ergodic theory, Lect.\ Notes in Math.\ \textbf{1132}, Springer, Berlin (1985), 556--588. 
  \bibitem{V2} D.\ Voiculescu, Addition of certain non-commutative random variables, J.\ Funct.\ Anal.\ \textbf{66} (1986), 323--346. 
  \bibitem{V-D-N} D.\ Voiculescu, K.J.\ Dykema and A.\ Nica, \textit{Free Random Variables}, CRM Monograph Series, AMS, 1992. 
\end{thebibliography}
\end{document}